\documentclass[11pt]{article}
\usepackage{amssymb}
   \csname @addtoreset\endcsname{equation}{section}
\newtheorem{theorem}{Theorem}%[section]
\newtheorem{lemma}{Lemma}%[section]
\newtheorem{proposition}{Proposition}%[section]
\newtheorem{definition}{Definition}%[section]
\newtheorem{remark}{Remark}%[section]
\def\e{\varepsilon}

\def\defi{\stackrel{{\scriptscriptstyle \Delta}}{=}}

\def\N{{\cal N}}

\def\o{\omega}
\def\O{\Omega}

\def\F{{\cal F}}
\def\w{\widehat}

\def\Re{{\rm Re\,}}

\def\R{{\bf R}}

\def\Z{{\cal Z}}

\def\L{L}

\def\g{\gamma}
\def\C{{\bf C}}

\def\X{{\cal X}}
\def\t{\theta}
\def\oo{\bar}

\def\p{\partial}
\def\G{\Gamma}
\def\GG{{\cal G}}
\def\U{{\cal U}}
\def\V{{\cal V}}

\def\M{{\cal M}}

\def\L{{\cal L}}

\newcommand{\be}{\begin{equation}}
\newcommand{\ee}{\end{equation}}
\newcommand{\bd}{\begin{displaymath}}
\newcommand{\ed}{\end{displaymath}}
\newcommand{\ba}{\begin{array}{ll}}
\newcommand{\ea}{\end{array}}
\newcommand{\baa}{\begin{eqnarray}}
\newcommand{\eaa}{\end{eqnarray}}
\newcommand{\baaa}{\begin{eqnarray*}}
\newcommand{\eaaa}{\end{eqnarray*}}
\font\sm=cmr10
\def\oo{\bar}

\def\K{{\cal K}}
\def\iomega{i\o}
\date{Revised version: May 7th, 2009; the first version: November 17, 2008 }
\title{
Predictability on finite horizon for processes with exponential
decrease
 of energy on higher frequencies %XXX
 }
\author{ Nikolai Dokuchaev\\ {\sm  Department of
Mathematics, Trent University, Ontario, Canada}}
 
\begin{document}
 \vspace{-0.5cm}
 \maketitle
\begin{abstract} The paper presents sufficient conditions of predictability
for continuous time processes in deterministic setting. We found
that processes with exponential decay on energy for higher
frequencies are predictable in some weak sense on some finite time
horizon defined by the rate of decay. Moreover, this predictability
can be achieved uniformly over classes of processes. Some explicit
formulas for predictors are suggested.
%\\
%PACS 2008 numbers: 02.30.Mv, %    Approximations and expansions
%02.30.Nw,  %  Fourier analysis
%02.30.Yy, %    Control theory
%07.05.Mh,  %  Neural networks, fuzzy logic, artificial intelligence
%07.05.Kf %Data analysis: algorithms and implementation; data management
\end{abstract}
\section{Introduction}
 We study pathwise predictability of
continuous time processes in deterministic setting and in the
framework of the frequency analysis. It is well known that certain
restrictions on frequency distribution can ensure additional
opportunities for prediction and interpolation of the processes. The
classical result is Nyquist-Shannon-Kotelnikov interpolation theorem
for the band-limited  processes.  There are related predictability
results;
 see, e.g., Wainstein and Zubakov (1962), Beutler (1966),
Brown(1969),
 Slepian (1978), Knab
(1981), Papoulis (1985), Marvasti (1986), Vaidyanathan (1987), Lyman
{\it et al} (2000, 2001). These works considered predictability of
single processes, and the crucial assumption was that the processes
are  band-limited; the predictors were non-robust with respect to
small noise in high frequencies; see, e.g., the discussion in
Chapter 17 from Higgins (1996).
\par
We study some special weak predictability of continuous time
processes. Instead of predictability of the original processes, we
study predictability of sets of anticausal convolution integrals for
a wide enough classes of kernels.  This version of predictability
was introduced in Dokuchaev (2008) for band-limited processes; it
allowed to establish uniform predictability in this weakened sense
over classes of  band-limited and high-frequency processes. In the
present paper, we established some predicability for continuous time
processes with exponential decay of energy on the higher
frequencies. It allows to consider processes that are not
band-limited. More precisely, we obtain a sufficient condition of
uniform weak predictability on prediction horizon $T$ over some
classes of processes with exponential decay of energy on higher
frequencies $\o\to \pm\infty$, when the energy is decreasing faster
than $e^{-T|\o|}$. An alternative formulation of this condition in
time domain is also given, The predictors are obtained explicitly in
the frequency domain via their transfer function. These predictors
are defined entirely by the kernel of the convolution integral and
their choice is independent from the characteristics of the
particular input processes.
\section{Problem setting and definitions}
Let $x(t)$ be a currently observable continuous time process,
$t\in\R$. The goal is to estimate, at a current time $t$, the values
$y(t)=\int_{t}^{t+T} k(t-s)x(s)ds$, where $k(\cdot)$ is a given
kernel, and $T>0$ is a given prediction horizon. At any time $t$,
the predictors use historical values of the observable process
$x(s)|_{s\le t}$.
\par
 We consider only linear predictors in the form $\w y(t)=\int_{-\infty}^t\w
k(t-s)x(s)ds$, where $\w k(\cdot)$ is a kernel that has to be found.
We will call $\w k$ a predictor or predicting kernel.
\par
Let us describe admissible classes of $k$ and $\w k$.
\par
 Let  $\R^+\defi[0,+\infty)$, $\C^+\defi\{z\in\C:\
\Re z> 0\}$, $i=\sqrt{-1}$.
\par
For $x\in  L_2(\R)\cup L_1(\R)$, we denote by $X=\F x$ the function
defined on $i\R$ as the Fourier transform of $x$; $$X(i\o)=(\F
x)(i\o)= \int_{-\infty}^{\infty}e^{-i\o t}x(t)dt,\quad \o\in\R.$$ If
$x\in L_2(\R)$, then $X$ is defined as an element of $L_2(\R)$ (more
precisely, $X(i\cdot)\in L_2(\R)$).
\par
For $x(\cdot)\in L_2(\R)$ such that $x(t)=0$ for $t<0$, we denote by
$\L x$  the Laplace transform \baa\label{Up} X(p)=(\L
x)(p)\defi\int_{0}^{\infty}e^{-p t}x(t)dt, \quad p\in\C^+. \eaa
\par
Let  $H^r$ be the Hardy space of holomorphic on $\C^+$ functions
$h(p)$ with finite norm
$\|h\|_{H^r}=\sup_{s>0}\|h(s+i\o)\|_{L_r(\R)}$, $r\in[1,+\infty]$
(see, e.g., Duren (1970)).
\begin{definition}
For $T>0$, we denote by $\K(T)$ the set of functions $k:\R\to\R$
such that $k(t)=0$ for $t\notin [-T,0]$ and such that $k\in
L_\infty(\R)$.
\end{definition}
\begin{definition}
Let $\w\K$ be the class of functions $\w k:\R\to\R$ such that $\w k
(t)=0$ for $t<0$ and such that $K(\cdot)=\L\w k\in H^2\cap
H^\infty$.
\end{definition}
\par
We consider below  $k\in\K(T)$ and $\w k\in\w\K$.
\begin{definition}\label{def1}
Let  $\oo\X$ be a class of processes $x(\cdot)$ from $L_2(\R)\cup
L_1(\R)$. Let $r\in[1,+\infty]$.
\begin{itemize}
\item[(i)]
 We say that the class $\oo\X$ is  $L_r$-predictable in
the weak sense with the prediction horizon $T$ if, for any
$k(\cdot)\in\K(T)$, there exists a sequence $\{\w
k_m(\cdot)\}_{m=1}^{+\infty}=\{\w
k_m(\cdot,\oo\X,k)\}_{m=1}^{+\infty}\subset \w\K$ such that $$
\|y-\w y_m\|_{L_r(\R)}\to 0\quad \hbox{as}\quad
m\to+\infty\quad\forall x\in\X, $$ where \baaa y(t)\defi
\int_t^{t+T}k(t-s)x(s)ds,\qquad \w y_m(t)\defi \int^t_{-\infty}\w
k_m(t-s)x(s)ds.\label{predict} \eaaa The process $\w y_m(t)$ is the
prediction of the process $y(t)$ which describes depends on the
future values of process $x(s)|_{s\in[t,t+T]}$.
\item[(ii)] Let the set $\F(\oo\X)\defi \{X=\F x,\quad x\in\oo \X\}$  be provided with a norm $\|\cdot\|$.
 We say that the class $\oo\X$ is  $L_r$-predictable in
the weak sense with the prediction horizon $T$ uniformly with
respect to the norm $\|\cdot\|$, if, for any $k(\cdot)\in\K(T)$,
there exists a sequence $\{\w k_m(\cdot)\} =\{\w k_m
(\cdot,\X,k,\|\cdot\|,\e)\}\subset \w\K$ such that $$ \|y-\w
y\,\|_{L_r(\R)}\to 0\quad\hbox{uniformly in} \quad \{x\in\oo\X:\
\|X\|\le 1\}. $$ Here  $y(\cdot)$  and $\w y_m(\cdot)$ are the
same as above.
\end{itemize}
\end{definition}
\section{The main result}
 For  $q\in\{1,2\}$, let $\X(q)=\X(q,T)$ be the  set of
processes $x(\cdot)\in L_2(\R)\cup L_1(\R)$ such that
\baaa\int_{-\infty}^{+\infty}e^{qT|\o|}|X(\iomega)|^q d\o <+\infty,
\quad X(\iomega)=\F x\eaaa
\par
 For $\O>0$, set $D(\O)\defi\R\backslash (-\O,\O)$.
\par
Clearly, if $x(\cdot)\in\X(q,T)$, then
\baaa\int_{D(\O)}e^{qT|\o|}|X(\iomega)|^q d\o\to
0\quad\hbox{as}\quad \O\to +\infty. \eaaa It can be seen also that,
for any $T>0$, the class $\X(q,T)$ includes all  band-limited
processes $x$ such that $X(i\o)=\F x\in L_q(\R)$, $q\in\{1,2\}$.
\begin{theorem}\label{ThM} Let
$q\in\{1,2\}$. Set $r=r(q)=q(q-1)^{-1}$ (in particular, $r=\infty$
for $q=1$ and $r=2$ for $q=2$). \begin{itemize}\item[(i)] The class
$\X(q,T)$ is $L_r$-predictable in the weak sense with the prediction
horizon $T$.
\item[(ii)] Let $\U(q)=\U(q,T)$   be a class of processes
$x(\cdot)\in \X(q,T)$ such that
 \baaa\int_{D(\O)}e^{qT|\o|}|X(\iomega)|^q
d\o\to 0\quad\hbox{as}\quad \O\to +\infty\quad\hbox{uniformly
on}\quad x(\cdot)\in \U(q). \eaaa Then this class $\U(q,T)$ is
$L_r$-predictable in the weak sense with the prediction horizon $T$
uniformly with respect to the norm $\|\cdot\|_{L_q(\R)}$.
\end{itemize}
\end{theorem}
\par
Some alternative descriptions and examples of sets $\U(q,T)$ are
given below.
\par  The question arises how to find the predicting kernels. In
the proof of Theorem \ref{ThM}, a possible choice of the kernels is
given explicitly in the frequency domain, i.e., via the transfer
functions.
\begin{remark}
In Dokuchaev (2008), similar weak predictability with infinite
horizon  was introduced and established for models where an ideal
low-pass filter exists; the predictors used in this paper were
different from the ones presented below. Theorem \ref{ThM} allows to
extend this weak predictability on the case when the filters are not
ideal but allow exponentially decay of energy on higher frequencies.
\end{remark}
\begin{remark}
The case when processes $k(\cdot)\in L_2(R)\backslash
L_{\infty}(\R)$ can be also covered. In this case, we have to
require that $x\in L_2(\R)$.
\end{remark}
%\begin{remark}\label{remE}
\subsection*{On possibility of extrapolation}
The weak predictability introduced in Theorem \ref{ThM} does not
ensure extrapolation of  the processes from $\U(q,T)$ in the
classical sense; even approximate extrapolation is not guaranteed.
Let us explain why.  Let $q=2$ and $T>0$ be given. Assume that the
values of the process $x(t)\in \U(2,T)$ are known for $t\le \tau$.
Consider a sequence of kernels $\{K_m\}_{m=1}^{+\infty}\subset
\K(T)$ that forms a orthonormal basis in the Hilbert space
$L_2(-T,0)$ (for instance, one may take Fourier series). For any
$\tau$, the function $x(t)|_{[\tau,\tau+T]}\in L_2(\tau,\tau+T)$ can
be represented as a series $x(t)=\sum_{m=1}^{\infty}f_mK_m(\tau
+T-t)$, where $f_m=\int_\tau^{\tau+T}K_m(\tau+T-s)x(s)ds$ are   the
corresponding Fourier coefficients (the series converges in
$L_2(\tau,\tau+T)$). Theorem \ref{ThM} ensures that, for any $\e>0$,
there exists a predictor such that the values of $f_m$
 can be predicted at time $t$ with
the error less or equal than $\e$ for all $m$. Unfortunately, it
does not help to predict the summa of infinite  series, even if this
$\e$ is small.
\subsubsection*{On  predicability and causality}
\label{SecClasses} It may appears that Theorem \ref{ThM} contradicts
to the obvious fact that a general process cannot be predicted in
any sense.   For instance, let $x_1(t)$ and $x_2(t)$ be two
processes such that $x_1(t)=x_2(t)=0$ for $t\le 0$ and such that
$x_1(t)\neq x_2(t)$ for $t>0$. Clearly, it is not possible to say
which process we observe at time $t=0$ using the values for $t\le
0$. (Some discussion and examples related to the predictability  can
be found in Chapter 17 from Higgins (1996)). However, it does not
contradict to Theorem \ref{ThM}. For instance, let
$x_k(t)=(-1)^kte^{-t}$ for $t>0$; it is easy to verify that the
Fourier transforms of the processes $x_k(\cdot)$  do not belong to
$\X(q)$, $q=1,2$. It reflects the lack of causality for these
process: the values for $t>0$ cannot be regarded as continuation of
some development started before $t=0$. In contrast, periodic,
almost-periodic, and band-limited processes have causality property
and therefore can be predicted.
\par
Theorem \ref{ThM} says that the predictability can be ensured for
some processes other than periodic, almost-periodic, or
band-limited. In particular, there is some causality for all
processes covered by this theorem, i.e., some signs of future
development are presented in the current time. Therefore, Theorem
\ref{ThM} can be interpreted as a new sufficient frequency condition
of causality for processes that are not periodic, almost-periodic,
or band-limited. In some cases, this causality  makes possible
predictability on some fixed finite horizon only. In particular,
these conditions are more restrictive for longer horizon.
\section{Sufficient conditions of predictability in time
domain}\label{secTD} In Theorem \ref{ThM}, conditions of
predictability are formulated in frequency domain. It can be useful
to add some sufficient conditions in time domain.
\par
 For $C>0$, consider a class $\M(C)$ of processes $x(t)\in
C^{\infty}(\R)$  such that there exists $M>0$ such that \baaa
&&\left\|\frac{d^kx}{dt^k}(\cdot)\right\|^2_{L_2(\R)}\le C^kM, \quad
k=0,2,4,6,....
\\
&&\frac{1}{2}\left(\left\|\frac{d^{k-1}x}{dt^{k-1}}(\cdot)\right\|^2_{L_2(\R)}
+\left\|\frac{d^{k+1}x}{dt^{k+1}}(\cdot)\right\|^2_{L_2(\R)}\right)\le
C^kM, \quad k=1,3,5,7,.... \qquad \forall x(\cdot)\in\M, \eaaa
\begin{proposition}\label{propD}
For any $T>0$ and $C>0$,  $\M(C)\subset \X(2,T)$. In particular, for
any $T>0$ and $C>0$, the class $\M(C)$ is
 predictable in the weak sense with the prediction horizon
$T$.
\end{proposition}
 \par
For $C>0$, consider a class $\N(C)$  of processes  $x(t)\in
C^{\infty}(\R)$  such that there exists $M>0$ such that \baaa
&&\hphantom{xxi}\left\|\frac{d^kx}{dt^k}(\cdot)\right\|^2_{L_2(\R)}\le
k!C^{-k}M, \quad k=0,2,4,6,....
\\
&&\frac{1}{2}\left(\left\|\frac{d^{k-1}x}{dt^{k-1}}(\cdot)\right\|^2_{L_2(\R)}
+\left\|\frac{d^{k+1}x}{dt^{k+1}}(\cdot)\right\|^2_{L_2(\R)}\right)\le
k!C^{-k}M, \quad k=1,3,5,7,.... \quad \forall x(\cdot)\in\N(C).
\eaaa
\begin{proposition}\label{propD!}
For any $C>2T$,  $\N(C)\subset \X(2,T)$. In particular,
 the class $\N(C)$ is
 predictable in the weak sense with the prediction horizon
$T<C/2$.
\end{proposition}
\section{Example: outputs of Gaussian filters} Let
 $\oo C>0$ and $\oo\o>0$ be given. Consider a class of processes
$\Z(\oo C,\oo\o)=\{z(\cdot)\}$ such that the Fourier transform
$Z(i\o)=\F z$ is defined in the class $L_2(\R)\cup L_{\infty}(\R)$
and $|X(\i\o)|\le\oo C$ for all $\o$ such that $|\o|>\oo\o$, where
$z\in\Z$.
\par
 Let
 $c_1>0$ and $v_1>0$ be given.
  Consider a set $\GG=\GG(c_1,v_2)$ of Gaussian filters with kernels
$k_\GG(t)=c\exp\left(-\frac{t^2}{v}\right)$ such that $v\ge v_1$ and
$|c|\le c_1$.
\par
Let $\V_\GG=\V_G(\oo C,\oo c,c_1,v_1)$ be the set of processes $x$
such that $x$ is a convolution of $z$ with a kernel $k_\GG$, where
$k_\GG\in \G(c_1,v_1)$, $z\in \Z(\oo C,\oo\o)$, i.e.,
$$ x(t)=\int_{-\infty}^{\infty}k_\GG(t-s)z(s)ds.
$$
 Note that generalized functions $v\in C(\R)^*$ are allowed to be
elements of $\Z$ (for instance, we include delta functions). In that
case, $\F z$ is still well defined, and  the corresponding process
$x=\F^{-1}X(i\o)$ is well defined in $L_2(\R)$, where
$X(i\o)=K_\GG(i\o)Z(i\o)$, $K_\GG(i\o)=\F k_\GG$.
\begin{proposition}\label{propG}  (i) $\V_\GG\subset \X(1)\cap\X(2)$; (ii) For $q=1,2$, for any $T>0$,
the class $\V_\GG$ is  predicted in the weak sense with the
prediction horizon $T$ uniformly with respect to the norm
$\|\cdot\|_{L_q(\R)}$.
\end{proposition}
\par
 For example, consider
processes \baa x(t)=\sum_{m=1}^N
c_{m}\exp\left(-\frac{(t-a_{m})^2}{v_{m}}\right).\label{gauss} \eaa
for some constants $N>0$, $c_{m}$, $a_{m}$, and $v_{m}>0$.
 By Proposition \ref{propG}, these processes belong to $\X(1)\cap\X(2)$,  and that any set of
these processes such that $$\{N\le C_1,\ |c_m|\le C_2,\ v_m\ge C_3,\
|a_m|\le C_4\}$$ forms a class $\U(q)$ with the properties required
in Theorem \ref{ThM}, for any given set of positive
$C_1,C_2,C_3,C_4$. Therefore, these processes can be predicted and
uniformly predicted in the weak sense of Theorem \ref{ThM}. In
particular, it is possible to predict for any $T>0$ the values
$\int_t^{t+T}k(t-s)x(s)ds$  using the values for $\t<t$ and the
predictors defined in the proof of Theorem \ref{ThM}. If $a_{m}>0$
and $v_m^{-1}$ are large enough, then the processes have sharp peaks
in $t>0$, and the values of $x(\t)$ are small for $\t\le 0$, and the
impact of the choice of $N$, $c_k$, and $a_k$ on
$y(\cdot)|_{(-\infty,0]}$ is small. However, this impact still
exists, ant it makes the prediction at time $t=0$ possible.
\subsubsection*{Special case: extrapolation of the still snapshot of the
temperature} We have regarded $x(t)$ as processes in time with the
time variable $t$. It is the most natural model for prediction.
However, there are other models where Theorem \ref{ThM} can be
applied. For instance, consider the problem of measurement of the
temperature on the one-dimensional rod. We consider the still
snapshot of the temperature rather than the dynamics of the process
of heat propagation. Let $x(t)$ be the temperature at the point with
the coordinate $t\in \R$. Let us assume that the temperature is
given as (\ref{gauss}) with $c_m>0$; this case corresponds to the
model when the heat was originated from $N$ point sources that were
applied at the points $t=a_m>0$ at past times defined by $v_m$. We
assume that $N,a_m,c_m$, and $v_m$ are unknown and non-observable.
Assume that the temperature can be measured in the points
$t\in(-\infty,0]$ only. The problem is to estimate integrals
$\int_0^{T}k(-s)x(s)ds$ using the observations at $t\in[0,+\infty)$
only, with $T>0$; in fact, it is a relaxed version of the
extrapolation problem. Theorem \ref{ThM} gives the solution, and the
"predictor" from the proof can be used.
 \section{Appendix: Proofs}
The proofs below  are very straightforward and do not use the
advanced theory of $H^p$-spaces; the existence of required
predictors is proved by presenting explicit transfer functions of
the predictors with desired properties.
\par
Let  $k(\cdot)\in\K(T)$ and $K(i\o)=\F k$. We assume here and below
that $\o\in\R$.
\par For
$\g\in\R$, $\g>0$,  set \baaa && g(p)\defi T\frac{\g-p}{\g+p}p,\quad
h(p)\defi g(p)-Tp,\quad V(p)\defi e^{h(p)}.
  \eaaa
\begin{lemma}\label{lemmaV}
\begin{itemize}
\item[(i)] $V(p)\in  H^{\infty}$ and
 $\w
K(i\o)=V(i\o)K(i\o)$ can be extended on $\C^+$ as function $\w
K(p)\in H^2\cap H^{\infty}$.
\item[(ii)]  $|V(i\o)|=\exp\left(\frac{2T\g
\o^2}{\g^2+\o^2}\right).$
\item[(iii)] $\sup_{\g>0}|V(i\o)|\le e^{T|\o|}.$
\item[(iv)] $V(i\o)\to 1$ as
$\g\to+\infty$ for all $\o\in\R$.
\item[(v)]
For any $\e>0$ and any $\O>0$, there exists $\g>0$ such that
$|V(i\o)- 1|\le\e$ for all $\o\in[-\O,\O]$.
\end{itemize}
\end{lemma}
\par
{\it Proof of Lemma \ref{lemmaV}}. Set $Q(i\o)=e^{-i\o T}K(i\o)$,
i.e., $K(i\o)=e^{i\o T}K(i\o)$. Clearly,
 \baaa Q(i\o)=e^{-i\o T} \int_{-T}^{0}e^{-i\o
t}k(t)dt=\int_{0}^{T}e^{-i\o (T-\t)}k(-\t)d\t
=\int_{0}^{T}e^{-i\o\tau}k(\tau-T)d\tau.\eaaa It follows that
$Q(i\o)$ can be extended on $\C^+$ as function $Q(p)\in H^2\cap
H^{\infty}$.
\par Further, $V(p)=e^{-Tp}e^{g(p)}$ and \baaa
g(p)=T\frac{\g-p}{\g+p}p=T\frac{-\g-p+2\g}{\g+p}p=-Tp+T\frac{2\g
p}{\g+p}. \eaaa It follows that $e^{g(p)}\in H^{\infty}$. Hence $\w
K(i\o)=V(i\o)K(i\o)=Q(i\o)e^{g(i\o)}$ can be extended on $\C^+$ as
function $\w K(p)\in H^2\cap H^{\infty}$.  Then statement (i)
follows.
\par
Further, \baaa &&g(i\o)=T
\frac{\g-i\o}{\g+i\o}i\o=T\frac{(\g-i\o)^2}{\g^2+\o^2}i\o=T\frac{\g^2-2\g
i\o-\o^2}{\g^2+\o^2}i\o. \eaaa
 Then
 \baaa
\Re h(i\o)= \Re g(i\o)= \frac{2T\g \o^2}{\g^2+\o^2}.  \label{re}
\eaaa
 Then statement (ii)  follows.
\par
Let us find the maximum of $\Re h(p)=\Re h(p,\g)$ in $\g\ge 0$. It
suffices to find $\g$ such that $\frac{\p}{\p \g}\Re h(i\o)=0$,
i.e., such that \baa \frac{\p}{\p \g}\left(\frac{2T\g
\o^2}{\g^2+\o^2}\right)=
\frac{2T\o^2(\g^2+\o^2)-4T\o^2\g^2}{(\g^2+\o^2)^2}=
2T\o^2\frac{\g^2+\o^2-2\g^2}{(\g^2+\o^2)^2}=0. \label{ree} \eaa It
is easy to see that (\ref{ree}) holds for $\g=|\o|$. For this
$\g=|\o|$, we have that $$\Re h(i\o)=\frac{2T\g
\o^2}{\g^2+\o^2}=\frac{2T |\o|\o^2}{2\o^2}=T |\o|.$$  Hence (iii)
follows.
\par
We have that \baaa h(p)=h(p,\g)=-Tp\frac{2p}{\g+p}. \eaaa Hence
$h(i\o,\g)\to 0$ as $\g\to +\infty$  for any $w\in\R$. Then
statement (iv) follows.
\par
Further, it follows from continuity of the exponent function that
there exists a function $\psi(\cdot):(0,+\infty)\to (0,+\infty)$
such that if $|h(i\o)|\le \psi(\e)$ then $|V(i\o)-1|<\e$. Let an
arbitrarily small $\e>0$ and an arbitrarily large $\O>0$ be given.
Take $\g=\g(\e)\ge 2T\O^2\psi(\e)^{-1}$, then \baaa
|h(i\o)|^2=\frac{4T^2\o^4}{\g(\e)^2+\o^2}\le \psi(\e)^2\quad
\forall \o\in[-\O,\O], \eaaa i.e., \baaa |V(i\o)-1|\le \e\quad
\forall \o\in[-\O,\O].\eaaa
  Then statement (v) follows.
 This completes the proof of Lemma \ref{lemmaV}. $\Box$
\par
{\it Proof of Theorem \ref{ThM}}. It suffices to present a set of
predicting kernels $\w k$ with desired properties. Let
$V(\cdot)=V(\g,\cdot)$ be as defined in Lemma \ref{lemmaV}. Set $\w
K(i\o)\defi V(i\o)K(i\o)$. Let the predicting kernels be defined as
$\w k(\cdot)=\w k(\cdot,\g(\e))=\F^{-1}\w K(i\o)$.
\par
 For  $x(\cdot)\in L_q(\R)$,
let
 $X(i\o)\defi \F x$,
 $Y(i\o)\defi \F y=K(i\o)X(i\o)$. Set $\w Y(i\o)\defi \w K(i\o)X(i\o)=V(i\o)Y(i\o)$ and $\w y=\F^{-1}\w Y$.
\par
Let us prove (i).  Since $K(i\o)\in L_{\infty}(\R)$ and $X(i\o)\in
L_q(\R)$, we have that $Y(i\o)=K(i\o)X(i\o)\in L_{q}(\R)$ and $\w
Y\in L_{q}(\R)$.  By Lemma \ref{lemmaV}(iv), it follows that  \baa\w
Y(i\o)\to Y(i\o)\quad\hbox{for a.e.}\quad \o\in \R
\quad\hbox{as}\quad\g\to +\infty. \label{YY}\eaa
 We have that $e^{T|\o|}X\in L_q(\R)$, $K(i\o)\in
L_{\infty}(\R)$ and \baa&&|\w K(i\o)-K(i\o)|\le
|V(i\o)-1||K(i\o)|\le 2e^{T|\o|} |K(i\o)|,\quad \o\in
D,\label{d1}\\ &&|\w Y(i\o)-Y(i\o)|\le
2e^{T|\o|}|K(i\o)||X(i\o)|,\quad \o\in D.\label{d2} \eaa
 By (\ref{YY}),(\ref{d2}), and by Lebesque Dominance Theorem, it follows that
\be\|\w Y-Y\|_{L_q(\R)}\to 0,\quad\hbox{i.e.,}\quad\|\w
y-y\|_{L_r(\R)}\to 0 \quad \o\in \R \quad\hbox{as}\quad\g\to
+\infty. \label{1s}\ee
\par
Let us prove (ii). Let $\e>0$ be given, and let
$\O(\e)\defi\e^{-1}$. By Lemma \ref{lemmaV}(v), there exists
$\g=\g(\e)>0$ such that $|V(i\o)-1|^q\le\e$ for all
$\o\in[-\O(\e),\O(\e)]$. For this $\g=\g(\e)$, we have \baa &&\|\w
Y(i\o)-Y(i\o)\|^q_{L_q(\R)}\nonumber\\&&\le
\int_{-\O(\e)}^{\O(\e)}|V(i\o)-1|^q|K(i\o)|^q|X(i\o)|^qd\o
+\int_{D(\O(\e))}|V(i\o)-1|^q|K(i\o)|^q|X(i\o)|^qd\o \nonumber\\
&&\le
\|K(i\o)\|^q_{L_\infty(\R)}\left(\e^q\|X(i\o)\|^q_{L_q(\R)}+2\int_{D(\O(\e))}e^{qT|\o|}|X(i\o)|^qd\o\right),
\label{4s}\eaa Take $\e\to 0$. By (\ref{4s}), it follows that
$\|\w Y(i\o)-Y(i\o)\|^q_{L_q(\R)}\to 0$ and $\|\w
y-y\|_{L_r(\R)}\to 0$ uniformly over $\U(q)\cap \{x(\cdot):\
\|X(i\o)\|_{L_q(\R)}\le 1\}$.
\par
By (\ref{1s}),(\ref{4s}), it follows  that the predicting kernels
$\w k(\cdot)=\w k(\cdot,\g(\e))=\F^{-1}\w K(i\o)$ are such as
required. This completes the proof of Theorem \ref{ThM}.
 $\Box$
 \par {\it Proof of
Proposition \ref{propD}}.  Let $x(\cdot)\in \M(C)$. Since
$\exp(2|\o|T)=\sum_{k=0}^{\infty}\frac{(2|\o|T)^k}{k!}$ and
$|\o|^k\le (|\o|^{k-1}+|\o|^{k+1})/2$, $k\ge 1$,  we have that \baa
\int_{-\infty}^{+\infty}e^{2|\o|T}|X(i\o)|^2d\o=
\sum_{k=0}^{\infty}\frac{(2T)^k}{k!}
\int_{-\infty}^{+\infty}|\o|^{k}|X(i\o)|^2d\o \le 2\pi
M\sum_{k=0}^{\infty}\frac{(2TC)^k}{k!}<+\infty.\label{chain1} \eaa
It follows that $x(\cdot)\in \X(2,T)$ for all $T>0$.  By Theorem
 \ref{ThM}(i), the required predictability holds.  $\Box$
\par
{\it Proof of Proposition \ref{propD!}}.  Let $x(\cdot)\in \N(C)$.
Similarly to (\ref{chain1}), we obtain that  \baaa
\int_{-\infty}^{+\infty}e^{2|\o|T}|X(i\o)|^2d\o=
\sum_{k=0}^{\infty}\frac{(2T)^k}{k!}
\int_{-\infty}^{+\infty}|\o|^{k}|X(i\o)|^2d\o \le 2\pi
M\sum_{k=0}^{\infty}\frac{(2T)^k}{k!}\frac{k!}{C^k}<+\infty. \eaaa
It follows that $x\in \X(2,T)$.  By Theorem
 \ref{ThM}(i), the required predictability holds.   $\Box$
 \par
{\it Proof of Proposition \ref{propG}}.  It is known that
$K_{\GG}(i\o)=c\sqrt{\pi v}\exp(-v\pi^2\o^2)$. It follows that
$e^{|\o|T}X(i\o)\in L_2(\R)\cap L_1(\R)$ and $x\in \X(1)\cap \X(2)$.
It follows also that condition (ii) in Theorem
 \ref{ThM} is satisfied for
class $\V_{\GG}$ for any $T>0$. $\Box$
%\section{Concluding remarks}
 \begin{remark}
The predictors introduced above  are stable, since the corresponding
transfer functions belong to $H^2\cap H^{\infty}$. In addition,
these predictors are robust with respect to the deviations of the
process that are small in the weighted norm generated by the
definition of the space $\X(q)$.
\end{remark}
\begin{remark}
Formally, the predictors described above require the past values
of $x(s)$ for all $s\in(-\infty,t]$, but it is not too
restrictive, since $\int_{-\infty}^t\w k(t-s)x(s)ds$ can be
approximated by $\int_{-M}^t\w k(t-s)x(s)ds$ for large enough
$M>0$. In addition, the corresponding transfer functions can be
approximated by rational fraction polynomials, and more general
kernels $k$ can be approximated by kernels from $\K(T)$.
\end{remark}

%\end{remark}
%\end{document}

\end{document}